\documentclass[11pt]{article}
\usepackage{mathrsfs}
\usepackage{latexsym,amsfonts,amssymb,amsthm}
\title{{\bf Graded modules for
Virasoro-like algebra}
\thanks{ Supported by the National Science Foundation of China
(No. 10371100, 10471096), the China Postdoctoral Science Foundation
(No. 20060390693), ``One Hundred Talents Program'' from University
of Science and Technology of China, and ``New Century Talents
Program" from Education Department of Fujian Province..
 $\mbox{$\dagger$
e-mail: linwq83@yahoo.com.cn}$ }}
\author{   Weiqiang\ Lin$^{1,2,\dagger}$\quad  Yucai\ Su$^1$\\
 {\small 1. Department of Mathematics, University of Science and Technology of China,}\\
{\small Hefei\ 230026, Anhui, China;}\\
{\small 2. Department of Mathematics, Zhangzhou
Teachers College,}\\
{\small Zhangzhou\ 363000, Fujian, China.}\\
}
\date{}
\begin{document}
\maketitle
 {\bf Abstract:} In this paper, we consider the classification of irreducible ${\bf
 Z}$-
 and  ${\bf Z}^2$-graded modules  with finite dimensional homogeneous subspaces over the
 Virasoro-like algebra.
 We first prove that such a module is a uniformly bounded module
 or a generalized highest weight module. Then we
 determine all generalized highest weight irreducible  modules. As
a consequence, we also determine all the modules with nonzero
center. Finally, we prove that there does not exist any nontrivial
${\bf Z}$-graded modules of intermediate series.



\begin{center}
\section*{$\S 1$\quad \bf Introduction}
\end{center}

The Virasoro algebra is playing an increasingly important role in
both mathematics and physics. A book on conformal field theory by Di
Francesco, Mathieu and Senechal [1] gives a great detail on the
connection between the Virasoro algebra and physics. It is
well-known that the Virasoro algebra acts on any  highest weight
module (except when the level is negative of dual coxeter number) of
the affine Lie algebra through the use of famous Sugawara operators.
Since  1992 when Mathieu [2] gave a classification of Harish-Chandra
modules of the Virasoro algebra,  it is natural for mathematicians
to generalize the theory of this Lie algebra. From an algebra point
of view, the Virasoro algebra can be regarded as the universal
central extension of the Lie algebra of derivations (denoted by
$\mbox{Der}A$) on the ring of Laurent polynomials $A={\bf
C}[t,t^{-1}]$. The Virasoro algebra has nontrivial positive-energy
unitary representations only if the center is nonzero and this is
one of the reasons why this algebra is more interesting than the
algebra $\mbox{Der}A$. So a natural generalization is the Lie
algebra $\mbox{Der}A_{\nu}$  of derivations on the Laurent
polynomials ring $A={\bf C}[t_1^{\pm 1},\cdots,t_{\nu}^{\pm 1}]$ in
commuting variables. $\mbox{Der}A_{\nu}$ is also known as Lie
algebra of the group of diffeomorphisms of $\nu$-dimensional torus.
Several attempts have been made by physicists to give a Fock space
representation of $\mbox{Der}A_{\nu}$ and its extension (see, e.g.,
[3]). Attempts failed to produce  interesting results probably due
to the lack of proper definition of ``normal ordering''. The first
surprising result is that $\mbox{Der}A_{\nu}$ with $\nu\geq2$ is
centrally closed ([4]). So people were forced to search other
algebras (similarly to $\mbox{Der}A_{\nu}$) which admit nontrivial
central extensions. There are two kinds of algebras being found. One
is the so-called higher rank Virasoro algebras introduced in [5].
See also the papers [6--8]. Recently, Su [9] and Lu, Zhao [10] have
completed the classification of Harish-Chandra modules over the
higher rank Virasoro algebras. But the results of classification
turn out to be  disappointed as the Harish-Chandra modules are
simply induced by
 modules of intermediate series, and thus the center of higher rank Virasoro
algebras must act as zero on these modules.

The other is the Virasoro-like algebra introduced in [11]. It is the
universal central extension of the skew derivation Lie algebra over
the Laurent polynomials ring in two commuting variables. It can be
generated by three elements and it contains many standard Heisenberg
Lie subalgebras. Some relations between the Virasoro-like algebra
$L$ and the generalized Clifford algebras were given in [12]. The
algebra $L$ has many properties similar to the Virasoro and
Heisenberg algebras. Many papers are devoted to the study of this
algebra. The derivation Lie algebra of the centerless Virasoro-like
algebra $\bar L$ and the automorphism group of this derivation Lie
algebra were studied in [13], while the structure of automorphism
group of $\bar L$ was considered in [14]. A large class of uniformly
bounded ${\bf Z}^2$-graded module over $\bar L$ was constructed in
[15], and a necessary condition for a nonzero level ${\bf
Z}^2$-graded irreducible module over  $L$ to be a module with finite
dimensional homogeneous subspaces was given in [16]. The paper [17]
 presented a large class of generalized highest weight
${\bf Z}^2$-graded irreducible modules over $\bar L$, the paper [18]
constructed a class of graded irreducible highest weight modules
over $L$ and classified the ${\bf Z}$-graded irreducible $L$-modules
with nonzero center and finite dimensional homogeneous subspaces,
while the paper [19] determined the structure of the Verma modules
over $L$. In the present paper, we deal with the classification of
irreducible graded $L$-modules with finite dimensional homogeneous
subspaces by using the results on the irreducible modules of
Heisenberg algebra obtained in [20], and the results about the ${\bf
Z}$-graded $L$-modules and the ${\bf Z}^2$-graded $L$-modules given
in [18] and [16].

The paper is arranged as follows. In section 2 we recall the
concepts of the Virasoro-like algebra and its graded modules. We
also collect some results about the irreducible modules of
Heisenberg algebra which is crucial in the study of the
classification of the graded irreducible $L$-modules. In section 3
we first prove that a ${\bf Z}$-graded $L$-module must be either a
uniformly bounded, or a generalized highest weight module. Then we
complete the classification of the irreducible generalized highest
weight ${\bf Z}$-graded modules with finite dimensional homogeneous
subspaces. In section 4 we first construct a class of irreducible
generalized highest weight ${\bf Z}^2$-graded modules with finite
dimensional homogeneous subspaces by using the results obtained in
the previous section. Then we prove that an irreducible generalized
highest weight ${\bf Z}^2$-graded $L$-modules must be isomorphic to
the modules constructed in the beginning of this section. Finally,
it is proved in Section 5 that there does not exist any nontrivial
${\bf Z}$-graded $L$-module of intermediate series.

\begin{center}
\item \section*{$\S 2$ \quad\bf The Virasoro-like algebra and its graded modules}
\end{center}

Throughout this paper we use ${\bf C},{\bf Z},{\bf Z_+},{\bf N}$ to
denote the sets of complex numbers, integers, nonnegative integers,
positive integers respectively. All spaces are over ${\bf C}$. Let
${\bf e}_1=(1,0),\; {\bf e}_2=(0,1)$ be the standard basis of ${\bf
C}^2$ and let $\Gamma={\bf Z}{\bf e}_1\oplus {\bf Z}{\bf e}_2$, an
additive group isomorphic to ${\bf Z}^2$. As usual, if
$u_1,\cdots,u_k$ are elements on some vector space, we use $\langle
u_1,\cdots,u_k\rangle $ to denote their linear span over $\bf C$.
Let $A_2={\bf C}[t_1^{\pm1},t_2^{\pm1}]$ be the Laurent polynomial
algebra with two variables. For ${\bf r}=r_1{\bf e}_1+r_2{\bf
e}_2\in\Gamma$, we denote $$t^{\bf r}=t_1^{r_1}t_2^{r_2}\in A_2.$$
We use $\frac{\partial}{\partial t_i}$ to denote the partial
derivative with respect to $t_i$ for $i=1,2$, and denote
$$d_i=t_i\frac{\partial}{\partial t_i}, \ i=1,2,\mbox{ \ \
and \ } D({\bf r})=t^{\bf
r}(r_2d_1-r_1d_2)=r_2t_1^{r_1+1}t_2^{r_2}\frac{\partial}{\partial
t_1}-r_1t_1^{r_1}t_2^{r_2+1}\frac{\partial}{\partial t_2}.$$ It is
clear that the derivation algebra of $A_2$ is $\mbox{Der}A_2=\langle
t^{\bf r}d_i\,|\,{\bf r}\in\Gamma,i=1,2\rangle $. One sees that the
subspace $$\overline{L}=\langle D({\bf r})\,|\,{\bf
r}\in\Gamma\rangle $$ forms a Lie subalgebra of $\mbox{Der}A_2$,
called the {\it skew derivation Lie algebra}, and also called the
{\it centerless Virasoro-like algebra} in [11]. For any ${\bf
m}=m_1{\bf e}_1+m_2{\bf e}_2,\ {\bf n}=n_1{\bf e}_1+n_2{\bf
e}_2\in\Gamma$, let $\mbox{det}({\bf m},{\bf n})=m_1n_2-m_2n_1$. One
can easily see that
$$[D({\bf m}),D({\bf n})]=
-\mbox{det}({\bf m},{\bf n})D({\bf m}+{\bf n}),
$$
and that $\overline{L}$ is a simple Lie algebra. The universal
center extension of $\overline{L}$ is the Lie algebra
$L=\overline{L}\oplus \langle c_1,c_2\rangle ,$ called the {\it
Virasoro-like algebra} (with center), with the following Lie
bracket:
$$
 [D({\bf m}),D({\bf n})]=
-\mbox{det}({\bf m},{\bf n})D({\bf m}+{\bf n}) +\delta_{{\bf m}+{\bf
n},{\bf 0}}h({\bf m}),$$
$$
c_1, c_2 \ \ \mbox{are central,}\eqno{(2.1)}
$$
for ${\bf m},{\bf n}\in\Gamma$, where$$h({\bf m})=m_1c_1+m_2c_2. $$
From the definition, one can directly deduces the following
result.\vspace{2mm}

{\it{\bf Lemma 2.1}\quad If ${\bf b}_1=b_{11}{\bf e}_1+b_{12}{\bf
e}_2,\ {\bf b}_2=b_{21}{\bf e}_1+b_{22}{\bf e}_2\in\Gamma$ then

{\rm(1)}\quad $L$ can be generated by $\{D({\bf b}_1),D({\bf
b}_2),D(-{\bf b}_1-{\bf b}_2)\}$ when $\{{\bf b}_1,{\bf b}_2\}$ is a
${\bf Z}$-basis of $\Gamma$.

{\rm(2)}\quad $\mathscr{H}=\langle D(k{\bf b}_2),h({\bf
b}_2)\,|\,k\in{\bf Z}\rangle $ is a standard Heisenberg subalgebra
of $L$.

{\rm(3)}\quad If ${\bf m}=m_1{\bf b}_1+m_2{\bf b}_2, {\bf n}=n_1{\bf
b}_1+n_2{\bf b}_2$, then
$$\begin{array}{ll}[D({\bf m}),D({\bf n})]= -\mbox{\rm det}({\bf b}_1,{\bf b}_2)\mbox{\rm
det}(m_1{\bf e}_1+m_2{\bf e}_2,n_1{\bf e}_1+n_2{\bf e}_2)D({\bf
m}+{\bf n})+\delta_{{\bf m}+{\bf n},{\bf 0}}h({\bf m}).
\end{array}$$}

Fix a ${\bf Z}$-basis ${\bf b}_1,{\bf b}_2$ of $\Gamma$. One can
regard $L$ as a {\it ${\bf Z}^2$-graded Lie algebra} by defining the
degree of the elements in $\langle D(m_1{\bf e}_1+m_2{\bf
e}_2)\rangle $ to be $(m_1,m_2)$ and the degree of the elements in
$\langle c_1,c_2\rangle $ to be $(0,0)$. One can also regard $L$ as
a {\it ${\bf Z}$-graded Lie algebra} by defining the $j\in{\bf Z}$
degree subspace of $L$ to be
 $$
 L_j=\langle D(j{\bf b}_1+m{\bf b}_2)\,|\,\;m\in{\bf Z}\rangle \oplus \delta_{j,0}\langle c_1,c_2\rangle.
\eqno{(2.2)}
$$

Next we recall the definitions of graded
$L$-modules. If an $L$-module $V=\oplus
 _{{\bf m}\in{\bf
Z}^2}V_{\bf m}$ satisfies
$$
D({\bf m})\cdot V_{\bf n}\subset V_{{\bf m}+{\bf n}},\quad\forall\
{\bf m},{\bf n}\in\Gamma, \eqno{(2.3)}
$$
then $V$ is called a {\it ${\bf Z}^2$-graded $L$-module} and $V_{\bf
m}$ is called a {\it homogeneous subspace of $V$ with degree} ${\bf
m}\in{\bf
Z}^2$. If $L$-module $V=\oplus
_{ m\in{\bf Z}}V_{m}$
satisfies
$$
L_m\cdot V_{n}\subset V_{ m+n},\quad\forall m,n\in{\bf Z},
\eqno{(2.4)}
$$
then $V$ is called a {\it ${\bf Z}$-graded $L$-module}, and $V_{m}$
is called a {\it homogeneous subspace of $V$ with degree} $m\in{\bf
Z}$.

A ${\bf Z}^2$- or ${\bf Z}$-graded module $V$ is called a {\it
quasi-finite graded module} if all homogeneous subspaces are finite
dimensional; a {\it uniformly bounded module} if there exists a
number $n\in{\bf N}$ such that all dimensions of the homogeneous
subspaces are $\le n$; a {\it module of the intermediate series} if
$n=1$. Denote the sets of irreducible quasi-finite ${\bf Z}^2$-  and
${\bf Z}$-graded $L$-modules by ${\cal O}_{{\bf Z}^2}$ and ${\cal
O}_{\bf Z}$ respectively. If a ${\bf Z}$-graded $L$-module $V$ is
generated by a vector $0\neq v\in V$ with $L_j\cdot v=0,\ \forall
j\in{\bf N}$, then $V$ is called a {\it highest weight module}.
Similarly, we have the notion of  a {\it lowest weight module}.
Furthermore, if there is a {\bf Z}-basis $B=\{{\bf b}'_1,{\bf
b}'_2\}$ of $\Gamma$, that is, ${\bf Zb}'_1+{\bf Zb}'_2=\Gamma$, and
a nonzero vector $v\in V_{n}$ (or $v\in V_{\bf n}$) such that
$D({\bf m})v=0,\forall\ {\bf m}\in{\bf Z}_{+}{\bf b}'_1+{\bf
Z}_{+}{\bf b}'_2$ (recall from definition that $D({\bf m})=0$ if
${\bf m}=(0,0)$), then $v$ is called a {\it generalized highest
weight vector} corresponding to the {\bf Z}-basis $B$, and $V$ is
called a {\it generalized highest weight module} with generalized
highest weight vector being of degree $n$ (or ${\bf n}$)
corresponding to $B$.

From the definition, one can easily see that the generalized highest
weight modules contain the highest weight modules and the lowest
weight modules as their special cases. Since the centers $c_1,\ c_2$
of $L$ must act as scalars  on any irreducible graded module, we
shall use the same symbols to denote these scalars.

 \vspace{3mm}Now we recall the construction and some basic properties of a class of modules over the
Heisenberg Lie algebra ${\mathscr H}$. Let $\psi$ be any linear
function on ${\mathscr H}$ such that $\psi(h({\bf b}_2))=0$. We
define the associative algebra homomorphism
$\overline{\psi}:U({\mathscr H})\rightarrow {\bf C}[t^{\pm 1}]$ such
that
$$
\overline{\psi}(D(k{ {\bf b}_2}))=\psi(D(k{ {\bf b}_2}))t^k,\;\ \
\overline{\psi}(h({\bf b}_2))=0,\;\ \ \forall\ k\in{\bf Z}^{*}:={\bf
Z}\backslash\{0\},
$$
where $U({\mathscr H})$ is the universal enveloping algebra of
${\mathscr H}$. Denote the image of $\overline{\psi}$ in ${\bf
C}[t^{\pm1}]$ by $A_{\overline{\psi}}$. Using $\overline{\psi}$, we
can define a ${\mathscr H}$-module structure on the space
$A_{\overline{\psi}}t^{i}$ for a given $i\in{\bf Z}$ as follows£º
$$
D(k{{\bf b}_2})t^m=\overline{\psi}(D(k{{\bf b}_2}))t^m,\;\ \ h({\bf
b}_2)t^m=0,\;\  \forall\ k\in {\bf Z}^{*},\ \ t^m\in
A_{\overline{\psi}}t^{i}.
$$

 \vspace{2mm} From the definition, we see  $A_{\overline{\psi}}\cdot t^{i}\simeq A_{\overline{\psi}}\cdot t^{j},\;
 \forall\ i,j\in{\bf Z} $ as ${\mathscr H}$-modules. Now, by Lemma 3.6 and Proposition 3.8 in [20],
 we have the following results.

\vspace{3mm}{ \it{\bf Theorem 2.2}\quad {\rm(1)}
$A_{\overline{\psi}} t^{i}$ is an irreducible ${\mathscr H}$-module
if and only if $A_{\overline{\psi}}t^{i}={\bf C}t^{i}$ $($denote
this module by $A_{0,i,0})$ or there exists $s\in {\bf Z}^*$ such
that $A_{\overline{\psi}} t^{i}={\bf C}[t^{\pm s}]t^{i}$ $($denote
this module by $A_{\overline{\psi},i,s})$.

{\rm(2)} If $V$ is an irreducible ${\bf Z}$-graded ${\mathscr
H}$-module with zero center, then there is a linear function $\psi$
over ${\mathscr H}$ and $i\in {\bf Z}$ such that $V\simeq
A_{\overline{\psi},i,s}$ or $V\simeq A_{0,i,0}$.

{\rm(3)} If $V$ is a uniformly bounded {\bf Z}-graded ${\mathscr
H}$-module, then the center $h({\bf b}_2)$ acts as zero on $V$.}

\vspace{2mm}The following lemma will be used in the next two
sections.\vspace{2mm}

{ \it{\bf Lemma 2.3}\quad An $L$-module $V$ is a generalized highest
weight modules if there is a  ${\bf Z}$-basis  ${\bf b}'_1, {\bf
b}'_2$ of $\Gamma$ and a homogeneous vector $v\neq 0$
 such that
 $D({\bf b}'_1)v=D({\bf b}'_2)v=0$.}

 {\bf Proof:~}
 By (2.1) and the assumption,  we can prove $D({\bf Nb}'_1+{\bf Nb}'_2)v=0$.
Thus for the ${\bf Z}$-basis ${\bf m}_1=3{\bf b}'_1+{\bf b}'_2,\;
{\bf m}_2=2{\bf b}'_1+{\bf b}'_2$ of $\Gamma$ we have $D({\bf
Z}_+{\bf m}_1+{\bf Z}_+{\bf m}_2)v=0$. Hence $V$ is a generalized
highest weight module by definition.\hfill$\Box$

 \begin{center}
\item \section*{$\S 3$ \quad\bf The ${\bf Z}$-graded modules over the Virasoro-like algebra $L$ }
\end{center}

{ \it{\bf Lemma 3.1}\quad An $L$-module $V$ is a generalized highest
weight module or a uniformly bounded module if it is a ${\bf
Z}$-graded module.}

{\bf Proof:}\quad Let $V=\oplus_{m\in{\bf Z}}V_m$. If $V$ is not a
generalized highest weight module, then for any $m\in{\bf Z},$ by
considering the following maps
$$
D(-m{\bf b}_1+{\bf b}_2):V_m\rightarrow V_0,\quad D((1-m){\bf
b}_1+{\bf b}_2):V_m\rightarrow V_1,
$$
we have
$$
\mbox{ker}D(-m{\bf b}_1+{\bf b}_2)\cap\mbox{ker}D((1-m){\bf
b}_1+{\bf b}_2)=0,
$$
by Lemma 2.3. Therefore $\mbox{dim}V_m\leq
\mbox{dim}V_0+\mbox{dim}V_1$. So $V$ is a uniformly bounded
module.\hfill$\Box$

\vspace{3mm}{ \it{\bf Lemma 3.2}\quad If $V$ is an irreducible
nontrivial generalized highest weight ${\bf Z}$-graded $L$-module
corresponding to the ${\bf Z}$-basis $B=\{{\bf b}'_1, {\bf b}'_2\}$
then

{\rm(1)}\quad For any $v\in V$ there is some $p\in{\bf N}$ such that
$D(m_1{\bf b}'_1+m_2{\bf b}'_2)v=0$ for all $m_1,m_2\geq p$.

{\rm(2)}\quad For any $0\neq v\in V$ and $m_1,m_2>0$, we have
$D(-m_1{\bf b}'_1-m_2{\bf b}'_2)v\neq 0$.}

{\bf Proof:}\quad We may assume that $v_0$ is a generalized highest
weight vector corresponding to a ${\bf Z}$-basis $B=\{{\bf b}'_1,
{\bf b}'_2\}$.

(1) By the irreducibility of $V$ and the PBW Theorem, there exists
$u\in U(L)$ such that \mbox{$v=u\cdot v_0$}, where  $u$ is a linear
combination of elements of the form $$u_n=D(i_1{\bf b}'_1+j_1{\bf
b}'_2)\cdots D(i_n{\bf b}'_1+j_n{\bf b}'_2).$$ Thus without loss of
generality, we may assume $u=u_n$. Take
$$\mbox{$
p_1=-\sum\limits_{i_s<0}i_s+1,\quad p_2=-\sum\limits_{j_s<0}j_s+1.
$}$$ By induction on $n$ one gets that $D(i{\bf b}'_1+j{\bf
b}'_2)\cdot v=0$ for $i\geq p_1$ and $j\geq p_2$, which gives the
result with $p=\mbox{max}\{p_1,p_2\}$.

(2) If there are $0\neq v\in V$ and $m_1,m_2>0$ with $D(-m_1{\bf
b}'_1-m_2{\bf b}'_2)v=0$, let $p$ be as in the proof of (1), then we
see $$D(-m_1{\bf b}'_1-m_2{\bf b}'_2),\;D({\bf b}'_1+p(m_1{\bf b}'_1
+m_2{\bf b}'_2)),\;D({\bf b}'_2+p(m_1{\bf b}'_1+m_2{\bf b}'_2)),$$
act trivially on $v$. By (2.1), the above elements generate the Lie
algebra $L$. So  $V$ is a trivial module, a contradiction with the
assumption.\hfill$\Box$

\vspace{3mm}{ \it{\bf Lemma 3.3}\quad An $L$-module $V\in{\cal
O}_{\bf Z}$ is a highest weight module or a lowest weight module if
it is a generalized highest weight module corresponding to a ${\bf
Z}$-basis $B=\{{\bf b}'_1, {\bf b}'_2\}$.}

{\bf Proof:}\quad For convenience, we suppose $V$ is a generalized
highest weight module with  generalized highest weight vector being
of degree $0$, corresponding to the ${\bf Z}$-basis ${\bf
b}'_1=b'_{11}{\bf b}_1+b'_{12}{\bf b}_2,$ ${\bf b}'_2=b'_{21}{\bf
b}_1+b'_{22}{\bf b}_2$. Let $a=b'_{11}+b'_{21}$ and denote
$$
\wp(V)=\{m\in {\bf Z}\mid V_m\neq 0\}.
$$
If necessary, by replacing ${\bf b}'_1,{\bf b}'_2$ by ${\bf
b}''_1=3{\bf b}'_1+{\bf b}'_2,\, {\bf b}''_2=2{\bf b}'_1+{\bf
b}'_2$, we can suppose $a\neq 0$.

First we prove that if $a>0$ then $V$ is a highest weight module.
 Let
$$
{\cal A}_i=\{j\in{\bf Z}\mid i+aj\in\wp(V)\},\;\; \forall\  0\leq
i<a,
$$
 then there is  $m_i\in{\bf Z}$ such that ${\cal
A}_i=\{j\in{\bf Z}\mid j\leq m_i\} $ or ${\cal A}_i={\bf Z}$ by
Lemma 3.2(2).

Set ${\bf b}={\bf b}'_1+{\bf b}'_2$. Now we prove ${\cal
A}_i\not={\bf Z}$ for all $0\leq i<a$. Otherwise, we may assume
${\cal A}_0={\bf Z}$. Thus we can choose $0\neq v_j\in V_{aj}$ for
any $j\in{\bf Z}$. By Lemma 3.2(1), we know that there is
 $p_{v_j}>0$ with
$$
D(s_1{\bf b}'_1+s_2{\bf b}'_2)\cdot v_j=0,\;\;\forall\
s_1,s_2>p_{v_j}. \eqno{(3.1)}
$$
Choose $\{k_j\in{\bf N}\mid j\in{\bf N}\}$ and $v_{k_j}\in V_{ak_j}$
such that
$$
k_{j+1}>k_j+p_{v_{k_j}}+2. \eqno{(3.2)}
$$
We prove that $\{D(-k_j{\bf b})\cdot v_{k_j}\mid j\in{\bf
N}\}\subset V_0$ is a set of linearly independent vectors, from this
we will get a contradiction which gives the result as required.
Indeed, for any $r\in{\bf N}$, there exists $a_r\in{\bf N}$ such
that $D(x{\bf b}+{\bf b}'_1)v_{k_r}=0,\;\forall x\geq a_r$ by Lemma
3.2(1). On the other hand, we have $D(x{\bf b}+{\bf b}'_1)\cdot
v_{k_r}\neq 0$ for any $x<-1$ by Lemma 3.2(2). Thus we can choose
$s_r\geq -2$ such that
$$
D(s_r{\bf b}+{\bf b}'_1)\cdot v_{k_r}\not=0,\quad D(x{\bf b}+{\bf
b}'_1)\cdot v_{k_r}=0,\;\forall x>s_r. \eqno{(3.3)}
$$
By (3.2) we have $k_r+s_r-k_j>p_{v_{k_j}}$ for all $1\leq j<r$.
Hence by (3.1) we know that for all $1\leq j<r$,
$$\begin{array}{ll}
D((k_r+s_r){\bf b}+{\bf b}'_1)\cdot D(-k_j{\bf
b})v_{k_j}&\!\!\!=[D((k_r+s_r){\bf b}+{\bf b}'_1),D(-k_j{\bf
b})]v_{k_j}
\\[5pt]&\!\!\!
=k_j\mbox{det}({\bf b}'_1,{\bf b}'_2)D((k_r+s_r-k_j) {\bf b}+{\bf
b}'_1)\cdot v_{k_j}=0. \end{array}$$ Now by (3.2) and (3.3), one
gets
$$\begin{array}{ll}
D((k_r+s_r){\bf b}+{\bf b}'_1)\cdot D(-k_r{\bf
b})v_{k_r}\!\!\!&=[D((k_r+s_r){\bf b}+{\bf b}'_1),D(-k_r{\bf
b})]v_{k_r}
\\[5pt]&
=k_r\mbox{det}({\bf b}'_1,{\bf b}'_2)D(s_r {\bf b}+{\bf b}'_1)\cdot
v_{k_r}\not=0. \end{array}$$ Hence if
$\sum_{j=1}^{r}\lambda_jD(-k_j{\bf b})\cdot v_{k_j}=0$ then
$\lambda_r=\lambda_{r-1} =\cdots=\lambda_{1}=0$ by the arbitrariness
of $r$. So we see $\{D(-k_j{\bf b})\cdot v_{k_j}\mid j\in{\bf
N}\}\subset V_0$ is a set of linearly independent vectors which
contradicts the fact that $V\in{\cal O}_{\bf Z}$. Therefore, for any
$0\leq i<a$, there is $m_i\in{\bf Z}$ such that
 ${\cal A}_i=\{j\in{\bf Z}\mid j\leq m_i\} $ which implies that
 $V$ is a highest weight module since
$\wp(V)=\bigcup_{i=0}^{a-1}{\cal A}_i$.

Similarly, one can prove that if $a<0$ then $V$ is a lowest weight
module.\hfill$\Box$

In order to complete the classification of the highest and lowest
weight ${\bf Z}$-graded $L$-module, we first give a triangular
decomposition of $L$ and construct a class of highest (lowest)
weight ${\bf Z}$-graded irreducible modules which are similar to
that described in [18]. Let
$$
{L}_{+}=\bigoplus\limits_{j>0}{ L}_{j},\;\;\;\;
{L}_{-}=\bigoplus\limits_{j<0}{L}_{j},
$$
where $L_j=\langle D(j{\bf b}_1+m{\bf b}_2)\,|\,\;m\in{\bf Z}\rangle
\oplus \delta_{j,0}\langle c_1,c_2\rangle $, then   ${L }$ has the
following triangular decomposition
$${ L}={ L }_{+}\oplus{ L }_0\oplus{ L }_{-}.$$
For any linear function $\psi$ over ${ L}_{0}$ with $\psi(h({\bf
b}_2))=0$, we define a one dimensional $({L}_{0}+{L}_{+})$-module
 ${\bf C}v_0$ as follows
$$
{L}_{j}v_0=0,\;\; xv_0=\psi(x)v_0,\;\;\forall\ j>0,\ x\in{L}_{0}.
$$
Then we get an induced ${L}$-module
$$
\overline{V}^{+}(\psi)=\mbox{Ind}^{{L}}_{{L}_{0}+{L}_{+}}{\bf C}v_0
=U({{L}})\otimes_{U({L}_{0}+{L}_{+})}{\bf C}v_0,
$$
where $U({L})$ is the universal enveloping algebra of ${L}$. Set the
degree of $v_0$ to be $0$ then $\overline{V}^{+}(\psi)$ becomes a
${\bf Z}$-graded module. It is obvious that
$\overline{V}^{+}(\psi)\simeq U({L}_{-})$ as vector spaces and
$\overline{V}^{+}(\psi)$ has an unique maximal proper submodule $J$.
Then we obtain an irreducible ${\bf Z}$-graded highest weight
${L}$-module
$$
\widehat{V}^{+}(\psi)=\overline{V}^+(\psi)/J.
$$
Similarly, we can define an irreducible lowest weight ${\bf
Z}$-graded ${L}$-module $\widehat{V}^{-}(\psi)$ for any linear
function $\psi$ over ${L}_{0}$ with $\psi(h({\bf b}_2))=0$.

\vspace{3mm}For convenience, we introduce the following definition,
which is similar to a definition in [17].\vspace{3mm}

{ \it{\bf Definition 3.4}\quad If $\psi$ is a nonzero linear
function over ${L}_0$ with $\psi(h({\bf b}_2))=0$ and there exist
$b_{10},b_{11},\cdots,b_{1s_1},\cdots, b_{r0},\cdots,b_{rs_r}\in{\bf
C},\ \ \alpha_1,\cdots,\alpha_r\in{\bf C}\setminus\{0\}$ such that
$$
\psi(D(k{\bf
b}_2))=\frac{(b_{10}+b_{11}k+\cdots+b_{1s_1}k^{s_1})\alpha_1^k+\cdots+
(b_{r0}+b_{r1}k+\cdots+b_{rs_r}k^{s_r})\alpha_r^k}{k},\;\; \forall\
k\in{\bf Z}^*,
$$
$$
\psi(h({\bf b}_1))=-\mbox{det}({\bf b}_1,{\bf
b}_2)(b_{10}+b_{20}+\cdots+b_{r0}),
$$
then $\psi$ is called an exp-polynomial function over ${L}_0$.
 }

\vspace{3mm}{ \it{\bf Remark 3.5}\quad Set $f_k=\psi(kD(k{\bf
b}_2)),\ f_0=-\mbox{det}({\bf b}_1,{\bf b}_2)\psi(h({\bf b}_1))$.
Then by a well-known
 combinatorial formula we see that
 $\psi\not\equiv 0$ is an exp-polynomial function over ${L}_0$ if and only if there are $a_0,\cdots,a_n\in{\bf C},$
 with
 $a_0a_n\neq 0$ such that
 $$
 \sum\limits_{i=0}^na_if_{k+i}=0,\;\psi(h({\bf b}_2))=0,\;\forall\ k\in{\bf
 Z}.
 $$
 }

By using a technique in [18], we prove the following theorem.

 \vspace{3mm}{ \it{\bf Theorem 3.6}\quad
The nontrivial ${\bf Z}$-graded $L$-module $\widehat{V}^{+}(\psi)$
is in ${\cal O}_{\bf Z}$ $($or
 $\widehat{V}^{-}(\psi)$ is in ${\cal O}_{\bf Z})$ if and only if $\psi$ is an exp-polynomial function over $L_0$.}

{\bf Proof:} For convenience, we introduce a linear map $\Psi: {\bf
C}[t_1^{\pm 1},t_2^{\pm 1}]\rightarrow L$ by defining
$\Psi(t_1^{m_1}t_2^{m_2})=D(m_1{\bf b}_1+m_2{\bf m}_2),\;\forall
m_1,m_2\in{\bf Z}$. By Remark 3.5, it is sufficient for us to prove
the following claim.

{\bf Claim 1.} $\widehat{V}^{+}(\psi)\in{\cal O}_{\bf Z}$ if and
only if there exists a nonzero polynomial
$P(t_2)=\sum_{i=0}^na_it_2^i\in{\bf C}[t_2]$ with $a_0a_n\neq 0$
such that
$$
\psi(\Psi(d_2(t_2^kP(t_2)))-a_{-k}\mbox{det}({\bf b}_1,{\bf
b}_2)h({\bf b}_1))=0,\;\;\forall k\in{\bf Z},
$$
where $a_k=0$ if $k\not\in\{0,1,\cdots,n\}$.

``$\Longrightarrow$''\quad Since $\dim V_{-1}<\infty$ and
$\Psi(t_1^{-1}t_2^i)\cdot v_0\in V_{-1}$ for all $i\in{\bf Z}$,
there exists a $k\in{\bf Z}$ and a nonzero polynomial
$P(t_2)=\sum_{i=0}^na_it_2^i\in{\bf C}[t_2]$ with $a_0a_n\neq 0$
such that
$$\Psi(t_1^{-1}t_2^kP(t_2))\cdot v_0=0.$$
Applying $\Psi(t_1t_2^s)$ for any $s\in{\bf Z}$ to the above
equality, we get
\begin{eqnarray*}
0&=&\Psi(t_1t_2^s)\cdot\Psi(t_1^{-1}t_2^kP(t_2))\cdot v_0=[D({\bf
b}_1+s{\bf b}_2), \sum\limits_{i=0}^na_iD(-{\bf b}_1+(k+i){\bf
b}_2)]\cdot v_0
\\&
=&\left(-\sum\limits_{i=0}^na_i\mbox{det}({\bf b}_1,{\bf
b}_2)(s+i+k)D((s+i+k){\bf b}_2)+a_{-s-k}h({\bf b}_1)\right)\cdot v_0
\\&
=&\left(\Psi(-\mbox{det}({\bf b}_1,{\bf
b}_2)d_2(t_2^{s+k}P(t_2)))+a_{-s-k}h({\bf b}_1)\right)\cdot v_0,
\end{eqnarray*}
which deduces the result as we hope.

``$\Longleftarrow$''\quad Since $L_-$ is generated by $L_{-1}$, and
$L_+$ is generated by $L_1$, one sees that
$$
L_{-1}\cdot V_{-i}=V_{-i-1},\;\; \forall i\in{\bf Z}_+,
$$
and $v=0$ if $L_1\cdot v=0$.

Next, we will show the following claim by induction on $l$.

{\bf Claim 2.} For any $l\in{\bf Z}_+$, there exists a nonzero
polynomial $P_l(t_2)=\sum_{i\in{\bf Z}}a_i^{(l)}t_2^i\in{\bf
C}[t_2]$ such that
$$
\left(\Psi(d_2(t_2^kP_l(t_2)))-a_{-k}^{(l)}\mbox{det}({\bf b}_1,{\bf
b}_2)h({\bf b}_1)\right)\cdot V_{-l}=0, \eqno{(3.4)}
$$$$
\Psi(t_1^{-1}t_2^kP_l(t_2))\cdot V_{-l}=0, \qquad\forall k\in{\bf
Z}. \eqno{(3.5)}
$$
By the assumption, (3.4) holds for $l=0$ with $P_0(t_2)=P(t_2)$. On
the other hand, we have
$$
\Psi(t_1t_2^s)\cdot \Psi(t_1^{-1}t_2^kP(t_2))\cdot
V_{0}=\left(\Psi(-\mbox{det}({\bf b}_1,{\bf
b}_2)d_2(t_2^{s+k}P(t_2)))+a_{-s-k}h({\bf b}_1)\right)\cdot
V_0=0,\quad \forall s,k\in{\bf Z},
$$
which deduces that (3.5) also holds for $l=0$ with
$P_0(t_2)=P(t_2)$. Thus the claim holds for $s=0$.

Suppose the claim holds for $l$. From (3.4) and (3.5), we have that
$$
\left(\Psi(d_2(Q(t_2)))-a_{Q}\mbox{det}({\bf b}_1,{\bf b}_2)h({\bf
b}_1)\right)\cdot V_{-l}=0, \eqno{(3.6)}
$$$$
\Psi(t_1^{-1}Q(t_2))\cdot V_{-l}=0, \eqno{(3.7)}
$$
for any $Q(t_2)\in{\bf C}[t^{\pm 1}]$ with $P_l(t_2)| Q(t_2)$, where
$a_Q$ is the constant in $Q(t_2)$.

Now let us consider the claim for $l+1$. Let
$P_{l+1}(t_2)=P_l(t_2)^3=\sum_{i\in{\bf Z}}a_i^{(l+1)}t_2^i$. Then
$P_l(t_2)|d_2(t_2^kP_{l+1}(t_2))$ and $P_l(t_2)|d_2\cdot
d_2(t_2^kP_{l+1}(t_2))$ for any $k\in{\bf Z}$. By using (3.6) and
(3.7), for any $s,k\in{\bf Z}$, we have that
\begin{eqnarray*}&\!\!\!\!\!\!\!\!\!\!\!\!\!\!\!\!\!\!\!\!&
\left(\Psi(d_2(t_2^kP_{l+1}(t_2)))-a_{-k}^{(l+1)}\mbox{det}({\bf
b}_1,{\bf b}_2)h({\bf b}_1)\right)\cdot (\Psi(t_1^{-1}t_2^s)\cdot
V_{-l}) \\&\!\!\!\!\!\!\!\!\!\!\!\!\!\!\!\!\!\!\!\!&=
[\Psi(d_2(t_2^kP_{l+1}(t_2))),\Psi(t_1^{-1}t_2^s)]\cdot V_{-l}
=-\mbox{det}({\bf b}_1,{\bf b}_2)\Psi(t_1^{-1}t_2^s(d_2\cdot
d_2\cdot (t_2^kP_{l+1}(t_2))))\cdot V_{-l}=0.
\end{eqnarray*}
Thus (3.4) holds for $l+1$. From the above equality and (3.7), we
deduce that, for any $m,k,s\in{\bf Z}$,
\begin{eqnarray*}&\!\!\!\!\!\!\!\!\!\!\!\!\!\!\!\!\!\!\!\!&
\Psi(t_1t_2^m).\Psi(t_1^{-1}t_2^kP_{l+1}(t_2)).(\Psi(t_1^{-1}t_2^s).V_{-l})
\\&\!\!\!\!\!\!\!\!\!\!\!\!\!\!\!\!\!\!\!\!&
=[\Psi(t_1t_2^m),\Psi(t_1^{-1}t_2^kP_{l+1}(t_2))].(\Psi(t_1^{-1}t_2^s).V_{-l})+
\Psi(t_1^{-1}t_2^kP_{l+1}(t_2)).\left(\Psi(t_1t_2^m).(\Psi(t_1^{-1}t_2^s).V_{-l})\right)
\\&\!\!\!\!\!\!\!\!\!\!\!\!\!\!\!\!\!\!\!\!&
=\left(\Psi(-\mbox{det}({\bf b}_1,{\bf
b}_2)d_2(t_2^{m+k}P_{l+1}(t_2)))+a_{-m-k}^{(l+1)}h({\bf
b}_1)\right).(\Psi(t_1^{-1}t_2^s).V_{-l})=0.
\end{eqnarray*}
So we have (3.4) also holds for $l+1$. This completes the proof of
Claim 2.

From Claim 2, we deduce that
$$\dim V_{-l-1}\leq \mbox{deg}P_{l+1}(t_2)\cdot\dim V_{-l},\quad
\forall l\in{\bf Z}_+.
$$
Thus $\widehat{V}^+(\psi)\in{\cal O}_{\bf Z}$. This completes the
proof of our theorem. \hfill$\Box$

 \vspace{3mm}{ \it{\bf Theorem 3.7}\quad If a ${\bf Z}$-graded $L$-module $V$ is in ${\cal O}_{\bf Z}$£¬
 then one and only one of the following cases holds.

 {\rm(1)}\quad $V$ is a uniformly bounded module.

 $(2)$\quad There exists an exp-polynomial function $\psi$ over $L_0$ such that $V\simeq
 \widehat{V}^{+}(\psi)$.

 $(3)$\quad There exists an exp-polynomial function $\psi$ over $L_0$ such that $V\simeq \widehat{V}^{-}(\psi)$.}

 {\bf Proof:}\quad If $V$ is not a uniformly bounded module then by Lemmas 3.1 and  3.3 we
obtain that $V$ is a highest weight module or a lowest weight
module. By the property  of
 a Verma module we know that there is a nonzero linear function $\psi$ over $L_0$ such that
 $V\simeq \widehat{V}^{+}(\psi)$ or
 $V\simeq \widehat{V}^{-}(\psi)$. Now by Theorem 3.6, case (2) or case (3) must hold.

Now we prove that only one of the above three cases holds. It is
obvious that only one of case (2) and
 case (3) holds. If case (2) holds then there is
 $k\in{\bf Z}^*$ with $\psi(D(k{\bf b}_2))\neq 0$ which implies that $D(k{\bf b}_2)v_0\neq 0$.
 Similar to the proof of part (2) in Lemma 3.2, one can easily
 deduce that $D(-{\bf b}_1)v_k\neq0$ for any nonzero homogeneous
 vector $v_k\in V_k$ of $V$. (In fact, if $D(-{\bf b}_1)v_k=0$ then
  $D((-k+1){\bf b}_1+{\bf b}_2), D((-k+1){\bf b}_1-{\bf
 b}_2)$ and $D(-{\bf b}_1)$ act trivially on $v_k$ by the
 construction of $\widehat{V}^+(\psi)$. Thus $V$ is a trivial module
 by the irreducibility of $V$ since from (2.1) one can easily check
 that the Lie algebra $L$ is generated by the above three elements.
 Hence we reach a contradiction.) Meanwhile
 $$
 {\cal B}=\{D(-{\bf b}_1)^jD((-n+j){\bf b}_1+k{\bf b}_2)v_0\,|\,0\leq
 j<n\}\subset \widehat{V}^+(\psi)_{-n},\ \forall n\in {\bf N}.
 $$
 Next we prove that ${\cal B}$ is a set of linear independent
 vectors in $\widehat{V}^+(\psi)_{-n}$. Indeed, if
 $$\mbox{$
 \sum\limits_{j=0}^{n-1}\lambda_jD(-{\bf b}_1)^jD((-n+j){\bf
 b}_1+k{\bf b}_2)v_0=0,
 $}$$
 then for any $0\leq i< n-1$ we have
 $$\begin{array}{ll}
 0\!\!\!&=D((n-i){\bf b}_1)\sum\limits_{j=0}^{n-1}\lambda_jD(-{\bf
 b}_1)^jD((-n+j){\bf b}_1+k{\bf b}_2)v_0\\[5pt]&
 =-\sum\limits_{j=0}^i\lambda_jk(n-i)\mbox{det}({\bf b}_1,{\bf
 b}_2)D(-{\bf b}_1)^jD((j-i){\bf b}_1+k{\bf b}_2)v_0,\end{array}
 $$
 which implies $\lambda_0=\cdots=\lambda_{n-1}=0$. Hence ${\cal B}$ is a set of linear independent vectors
 in $\widehat{V}^{+}(\psi)_{-n}$ and thus
 $$
 \mbox{dim}(\widehat{V}^{+}(\psi)_{-n})\geq n.
 $$
Therefore $\widehat{V}^{+}(\psi)$ is not a uniformly bounded module
by the arbitrariness of $n$. So only one of case (2) and (1) holds.

 Similarly we can prove that only one of case (1) and case (3) holds.\hfill$\Box$

 \vspace{3mm}{\bf Remark 3.8}\quad In Section 5, we will prove
 that there does not exists a nontrivial  ${\bf
 Z}$-graded $L$-module of the intermediate series. We conjecture that there does not exist any nontrivial
 uniformly bounded ${\bf Z}$-graded irreducible $L$-module.
\vspace{3mm}

 From Theorem 3.7, we can recover  Theorem 3.2 in [18] below.

 \vspace{3mm}{ \it{\bf Corollary 3.9}\quad If a
  center of $L$ acts as a nonzero scalar on the
${\bf Z}$-graded module $V$, then $V\in{\cal O}_{{\bf Z}}$ if and
only if there is an exp-polynomial function $\psi$ over $L_0$ such
that $V\simeq \widehat{V}^{+}(\psi)$ or  $V\simeq
\widehat{V}^{-}(\psi)$.}

 {\bf Proof:}\quad Set $V=\oplus_{n\in{\bf Z}}V_n$. We have $h({\bf b}_2)=0$ by considering
 $V_0$ as a finite dimensional $L_0$-module. Thus $h({\bf b}_1)\neq 0$ since a center of $L$ acts as a
 nonzero scalar on $V$ and $\{{\bf b}_1,{\bf b}_2\}$ is a ${\bf Z}$-basis of
 $\Gamma$. Therefore, by Theorem 2.2(3), the module $V$ is not a uniformly bounded module
 since the Lie subalgebra  $\langle D(m{\bf b}_1),h({\bf b}_1)\mid m\in{\bf Z}\rangle $ of $L$
 is isomorphic to ${\mathscr H}$ and its center $h({\bf b}_1)$ acts as a nonzero scalar on
 $V$. Thus we obtain the result as we
 wish by Theorem 3.7.\hfill$\Box$

 \begin{center}
\item \section*{$\S 4$ \quad\bf The ${\bf
Z}^2$-graded modules over the Virasoro-like algebra $L$ }
\end{center}

We first construct a class of irreducible generalized highest weight
${\bf Z}^2$-graded $L$-modules by using the ${\bf Z}$-graded
$L$-module $\widehat{V}^{+}(\psi)$. For any linear function $\psi$
over $L_0$ with $\psi(h({\bf b}_2))=0$, we set
 $V(\psi)=\widehat{V}^+(\psi)\otimes{\bf C}[t^{\pm1}]$, and define the actions of the elements of $ L $ on $V(\psi)$
 as follows
$$
D(i{\bf {\bf b}_1}+j{\bf {\bf b}_2})(v\otimes t^k)=(D(i{\bf {\bf
b}_1}+j{\bf {\bf b}_2})v)\otimes t^{k+j}, \eqno{(4.1)}
$$
where $ v\in\widehat{V}^+(\psi),\ i,j,k\in{\bf Z}$. For any
homogeneous element $v$ with degree $m$ in ${\bf Z}$-graded module
$\widehat{V}^{+}(\psi)$, define the degree of $v\otimes t^k$ to be
$m{\bf b}_1+k{\bf b}_2$. Then one can easily see that $V(\psi)$
becomes a ${\bf Z}^2$-graded $L$-module. Denote $U({\mathscr
H})\cdot(v_0\otimes t^i)$ by $W_i$, then $W_i\simeq
A_{\overline{\psi}\cdot t^i}$ as ${\mathscr H}$-modules. If $W_i$ is
an irreducible ${\mathscr H}$-module then $W_i\simeq
A_{\overline{\psi},i,s}$ or $W_i\simeq A_{0,i,0}$ by Theorem 2.2.
Thus by the construction of $L$-module $V(\psi)$ we know  that there
exists a unique maximal proper submodule either
$W_{\overline{\psi},i,s}$ or $W_{0,i,0}$ which intersects trivially
with $W_i$. Then we have the irreducible ${\bf Z}^2$-graded
$L$-module either
 $$
 \overline{V}(\psi,i,s)=V(\psi)/W_{\overline{\psi},i,s},\quad \mbox{or \ \ }\overline{V}(0,i,0)=V(\psi)/W_{0,i,0}.
 \eqno{(4.2)}
 $$
One can easily check that there exists $s\in{\bf N}$ such that
 $$
 W_i=v_0\otimes{\bf C}[t^{\pm s}]\cdot t^i\mbox{  if $\psi$ is an exp-polynomial function over $L_0$.
}
 $$
Hence by Theorem 2.2 we know  $W_i$ is an irreducible ${\mathscr
H}$-module in this instance. Now by Theorem 3.6 and the construction
of $\overline{V}(\psi,i,s)$ we obtain the following result.

 \vspace{3mm}{ \it{\bf Lemma 4.1}\quad If $\psi$ is an exp-polynomial function over $
 L_0$ then

 $(1)$\quad $\overline{V}(\psi,i,s)\in{\cal O}_{{\bf Z}^2}$;
 \qquad$(2)$\quad $V(\psi)\simeq \overline{V}(\psi,i,s)$
 if and only if $s=1$. }

 \vspace{3mm}{ \it{\bf Remark 4.2}\quad   $V(\psi)$ is not an irreducible $L$-module in
 general.}
 For example, we can define a linear function  $\psi$ over $L_0$
 such that
 $$\psi(D(k{\bf b}_2))=\frac{(-1)^k+1}{k},\;\forall\; k\in{\bf Z}^*,\quad
 \psi(h({\bf b}_1))=-2\det({\bf b}_1,{\bf b}_2),\;\psi(h({\bf b}_2))=0.
 $$
Then one can easily see $\psi$ is an exp-polynomial function over
$L_0$. On the other hand, by the PBW Theorem and the construction of
$L$-module $V(\psi)$, we can deduce
 $$
 W_0=v_0\otimes{\bf C}[t^{\pm 2}],\quad (U(L).(v_0\otimes t))\bigcap
 W_0=0.
 $$
Thus $U(L)\cdot(v_0\otimes t)$ is a nonzero proper submodule of
$V(\psi)$ which implies that $V(\psi)$ is not  irreducible.

 \vspace{3mm}{ \it{\bf Lemma 4.3}\quad A  ${\bf Z}^2$-graded $L$-module is
  a generalized highest weight module or a uniformly bounded module.}

 {\bf Proof:}\quad Suppose $V=\oplus_{{\bf m}\in{\bf Z}^2}V_{\bf m}$.
If $V$ is not a generalized highest weight module, then for any
$(m_1,m_2)\in{\bf Z}^2$,  by considering the linear maps
$$\begin{array}{ll}
D(-m_1{\bf e}_1+{\bf e}_2):&V_{(m_1,m_2)}\rightarrow V_{(0,m_2+1)},
\\[5pt]
D((1-m_1){\bf e}_1+{\bf e}_2):&V_{(m_1,m_2)}\rightarrow
V_{(1,m_2+1)},
\end{array} $$ we obtain by Lemma 2.3
$$
\mbox{ker}D(-m_1{\bf e}_1+{\bf e}_2)\cap\mbox{ker}D((1-m_1){\bf
e}_1+{\bf e}_2)=0.
$$
Thus
$$
\mbox{dim}V_{(m_1,m_2)}\leq
\mbox{dim}V_{(0,m_2+1)}+\mbox{dim}V_{(1,m_2+1)}. $$
 Now consider
the following linear maps
$$\begin{array}{ll}
D(-{\bf e}_1-m_2{\bf e}_2):&V_{(0,m_2+1)}\rightarrow V_{(-1,1)},
\\[5pt]
 D(-{\bf e}_1+(1-m_2){\bf e}_2):&V_{(0,m_2+1)}\rightarrow V_{(-1,2)}.
 \end{array}$$
By Lemma 2.3 we have
$$
\mbox{ker}D(-{\bf e}_1-m_2{\bf e}_2)\cap\mbox{ker}D(-{\bf
e}_1+(1-m_2){\bf e}_2)=0.
$$
Thus $\mbox{dim}V_{(0,m_2+1)}\leq
\mbox{dim}V_{(-1,1)}+\mbox{dim}V_{(-1,2)}$.

Similarly we can deduce
 $$\mbox{dim}V_{(1,m_2+1)}\leq
 \mbox{dim}V_{(0,1)}+\mbox{dim}V_{(0,2)}.$$
 Therefore $V$ is a uniformly bounded module.\hfill$\Box$

Now we recall the following proposition from [16].

\vspace{3mm}{ \it{\bf Theorem 4.4} $([16]$ Proposition $3.9)$\quad
If $V$ is a generalized highest weight ${\bf Z}^2$-graded $L$-module
and $V\in{\cal O}_{{\bf Z}^2}$, then there exists a ${\bf Z}$-basis
$B=\{{\bf b}_1,{\bf b}_2\}$ of $\Gamma$ and a linear function $\psi$
over $L_0=\langle m{\bf b}_2)\,|\,\;m\in{\bf Z}\rangle \oplus
\langle c_1,c_2\rangle $ with $\psi(h({\bf b}_2))=0$ such that $V$
is isomorphic to $\overline{V}(\psi,i,s)$ or $\overline{V}(0,i,0)$,
where $i\in{\bf Z}, s\in{\bf N}$ .}

 \vspace{3mm} { \it{\bf Lemma 4.5}\quad $(1)$\quad If $\psi(h({\bf b}_1))=0$ then
 $\overline{V}(0,i,0)$ is a trivial module.

$(2)$\quad If $\psi(h({\bf b}_1))\neq 0$ then
$\overline{V}(0,i,0)\not\in {\cal O}_{{\bf Z}^2}$.}

  {\bf Proof:}\quad (1) By the construction of $\overline{V}(0,i,0)$, one can easily deduce the result.

  (2) By definition we have
  $$\psi(D(m{\bf b}_2))=\psi(h({\bf b}_2))=0, \ \forall\ m\in{\bf Z},$$
  and
  $$
  {\cal B}_n=\{D(-{\bf b}_1+m{\bf b}_2)D(-{\bf b}_1-m{\bf b}_2)\cdot
  (v_0\otimes t^i)\mid 1\leq m\leq n\}\subset \overline{V}_{-2{\bf b}_1+i{\bf b}_2},
  \;\forall\ n\in{\bf N}.
  $$
 Since
 $$
 D({\bf b}_1+m{\bf b}_2)D(-{\bf b}_1-m{\bf b}_2)\cdot
 (v_0\otimes t^i)=h({\bf b}_1+m{\bf b}_2)\cdot(v_0\otimes t^i)\neq 0, \;\forall
 m\in{\bf Z},
 $$
  we have
  $$D(-{\bf b}_1-m{\bf b}_2)\cdot(v_0\otimes t^i)\neq 0.$$

Now we  prove that ${\cal B}_n$ is a set of linear independent
elements in $\overline{V}_{-2{\bf b}_1+i{\bf b}_2}$. If there are
$\lambda_1,\cdots,\lambda_n\in{\bf C}$ such that
  $$\mbox{$
  \sum\limits_{j=1}^n\lambda_j D(-{\bf b}_1+j{\bf b}_2)D(-{\bf b}_1-j{\bf b}_2)\cdot(v_0\otimes t^i)=0,
  $}$$
then for any $1\leq k\leq n$ we have
  $$\begin{array}{lll}
  0\!\!\!\!&=D({\bf b}_1-k{\bf b}_2)\sum\limits_{j=1}^n\lambda_j D(-{\bf b}_1+j{\bf b}_2)D(-{\bf b}_1-j{\bf b}_2)\cdot(v_0\otimes t^i)
  \\[5pt]&
  =\lambda_k \psi(h({\bf b}_1))D(-{\bf b}_1-k{\bf b}_2)\cdot(v_0\otimes t^i)+
  \sum\limits_{j=1}^n\lambda_j (j-k)^2(\mbox{det}({\bf b}_1,{\bf b}_2))^2D(-{\bf b}_1-k{\bf b}_2)\cdot(v_0\otimes
  t^i)
  \\[5pt]&
  =\lambda_k \psi(h({\bf b}_1))D(-{\bf b}_1-k{\bf b}_2)\cdot(v_0\otimes t^i)+
  \sum\limits_{j=1}^n\lambda_j (j-k)^2D(-{\bf b}_1-k{\bf b}_2)\cdot(v_0\otimes
  t^i).\end{array}
  $$
 On the other hand, for any $k>n$, we have
  $$\begin{array}{ll}
  0\!\!\!\!&=D({\bf b}_1-k{\bf b}_2)\sum\limits_{j=1}^n\lambda_j D(-{\bf b}_1+j{\bf b}_2)D(-{\bf b}_1-j{\bf b}_2)\cdot(v_0\otimes t^i)
  \\[5pt]&
  =\sum\limits_{j=1}^n\lambda_j (j-k)^2D(-{\bf b}_1-k{\bf b}_2)\cdot(v_0\otimes
  t^i).\end{array}
  $$
  Thus
  $$\mbox{$\lambda_k \psi(h({\bf b}_1))+
  \sum\limits_{j=1}^n$}\lambda_j (j-k)^2=0,\;\;\forall \ 1\leq k\leq n;
  \eqno{(4.3)}
  $$$$\mbox{$
  \sum\limits_{j=1}^n$}\lambda_j (j-k)^2=0,\;\;\forall \  k>n.
  \eqno{(4.4)}
  $$
 Hence (4.4) implies that the polynomial $
  \sum_{j=1}^n\lambda_j (j-x)^2=0$ has infinite many roots, and
  therefore $\sum_{j=1}^{n}\lambda_j(j-x)^2$ is a zero polynomial.
  This, together with (4.3), gives $\lambda_k=0,\ \forall 1\leq
  k\leq n$. Hence ${\cal B}_{n}$ is a set of linear independent
  vectors in $\overline{V}_{-2{\bf b}_1+i{\bf b}_2}$ which implies $\dim
  \overline{V}_{-2{\bf b}_1+i{\bf b}_1}\geq n,\ \forall n\in{\bf N}$. Thus
  $\overline{V}(0,i,0)\not\in{\cal O}_{{\bf Z}^2}$.
 \hfill$\Box$

 \vspace{3mm}{ \it{\bf Theorem 4.6}\quad If $V$ is a nontrivial generalized highest weight  ${\bf
 Z}^2$-graded $L$-module, then ${V}\in{\cal O}_{{\bf Z}^2}$ if and only if
there exists
 a ${\bf Z}$-basis $B=\{{\bf b}_1,{\bf b}_2\}$ of $\Gamma$ and an exp-polynomial function $\psi$ over
 ${L}_0=\langle m{\bf b}_2)\,|\,\;m\in{\bf Z}\rangle \oplus \langle c_1,c_2\rangle $ such that $V$ is
 isomorphism to $\overline{V}(\psi,i,s)$, where $i\in{\bf Z},\;s\in{\bf N}$.}

 {\bf Proof:}\quad The sufficiency follows directly from part (1) of Lemma 4.1. Now we prove the necessity.
 By Theorem 4.4 and Lemma 4.5 we know  there exist  a ${\bf Z}$-basis $B=\{{\bf b}_1,{\bf b}_2\}$ of $\Gamma$
and a linear function $\psi$ over $L_0$ with $\psi(h({\bf b}_2))=0$
such
 that $V$ is isomorphism to  $\overline{V}(\psi,i,s)$, where  $i\in{\bf Z}, s\in{\bf N}$.

 If $s=1$ then by Lemma 4.1(2) we have
 $\overline{V}(\psi,i,s)\simeq V(\psi)$. Thus by the construction of $V(\psi)$ we see that the
 dimensions of the homogeneous
 subspaces $\overline{V}(\psi,i,s)_{m{\bf b}_1+n{\bf b}_2}$ of $\overline{V}(\psi,i,s)$ equal to those of the
 corresponding homogeneous subspaces  $\widehat{V}^+(\psi)_m$ of $\widehat{V}^+(\psi)$.
 Then the result for the case $s=1$ follows directly from Theorem
 3.6.

 If $s\not=1$ then $\sum_{i=1}^{s}\overline{V}(\psi,i,s)$ is a quasi-finite ${\bf Z}^2$-graded module.
Moreover, since
 $$v_j=D(-{\bf b}_1+j{\bf b}_2)\cdot(v_0\otimes t^{-j})\in\Big(\mbox{$\sum\limits_{i=1}^{s}$}\overline{V}(\psi,i,s)\Big)_{-{\bf b}_1},\;\;\forall\ j\in{\bf Z},
 $$
 there exist an integer $k\in{\bf Z}$, and numbers $a_0,\cdots,a_n\in{\bf C}$ with $a_0a_n\neq 0$
 such that
 $ \sum_{j=0}^{n}a_jv_{k+j}=0$.  For convenience, we set
 $a_i=0$ for $i<0 \ \mbox{ and}\ i>n$. Then for any $s\in{\bf Z}$ we have
 $$\begin{array}{ll}
0\!\!\!\!&= D({\bf b}_1-s{\bf b}_2)\sum\limits_{j=0}^{n}a_jv_{k+j}
\\[5pt]&
=\sum\limits_{j=0}^{n}a_j(s-(k+j))\mbox{det}({\bf b}_1,{\bf b}_2)
\psi(D((k+j-s){\bf b}_2))(v_0\otimes t^{-s})+a_{s-k}\psi(h({\bf
b}_1))(v_0\otimes t^{-s}).\end{array}
$$
Thus
$$\mbox{$
\sum\limits_{j=0}^{n}$}a_j\psi((k+j)D((k+j){\bf
b}_2))-a_{-k}\mbox{det}({\bf b}_1,{\bf b}_2)\psi(h({\bf
b}_1))=0,\;\forall\ k\in{\bf Z}.
$$
 Denote
 $$
 f_k=\psi(kD(k{\bf b}_2)),\;\forall\ 0\neq k\in{\bf Z},\;f_0=-\mbox{det}({\bf
b}_1,{\bf b}_2)\psi(h({\bf b}_1)),
$$
 then
 $$
 \mbox{$\sum_{j=0}^{n}$}a_jf_{k+j}=0,\;\forall\ k\in{\bf Z}.
 $$
 Therefore by
 Remark 3.5 we see $\psi$ is an exp-polynomial function over $L_0$.\hfill$\Box$

 \vspace{3mm}By  similar arguments as in the proof of Theorem 3.7 and Corollary 3.9, one can prove the following
 two results by applying Lemma 4.3 and Theorems 4.6 and  2.2.

\vspace{3mm} { \it{\bf Theorem 4.7}\quad If ${\bf Z}^2$-graded
 $L$-module $V\in{\cal O}_{{\bf Z}^2}$, then one and only one of the following cases holds.

 $(1)$\quad $V$ is a uniformly bounded module.

 $(2)$\quad There exist
 a ${\bf Z}$-basis $B=\{{\bf b}_1,{\bf b}_2\}$ of $\Gamma$ and an exp-polynomial function $\psi$ over
 ${L}_0=\langle m{\bf b}_2)\,|\,\;m\in{\bf Z}\rangle \oplus \langle c_1,c_2\rangle $ such that $V$ is
 isomorphism to $\overline{V}(\psi,i,s)$, where $i\in{\bf Z},\;s\in{\bf N}$.}

\vspace{3mm}{ \it{\bf Corollary 4.8}\quad If $V$ is a  ${\bf
Z}^2$-graded $L$-module with nonzero centers, then $V\in{\cal
O}_{{\bf Z}^2}$ if and only if there exist
 a ${\bf Z}$-basis $B=\{{\bf b}_1,{\bf b}_2\}$ of $\Gamma$ and an exp-polynomial function $\psi$ over
 ${L}_0=\langle m{\bf b}_2)\,|\,\;m\in{\bf Z}\rangle \oplus \langle c_1,c_2\rangle $ such that $V$ is
 isomorphism to $\overline{V}(\psi,i,s)$, where $i\in{\bf Z},\;s\in{\bf N}$.}
\vspace{3mm}

 We would like to conclude this section by recalling
some results from [15] which show that there are many uniformly
bounded irreducible ${\bf Z}^2$-graded $L$-modules. A question one
may ask is: whether or not the modules constructed in the following
exhaust all the uniformly bounded irreducible ${\bf Z}^2$-graded
$L$-modules.

Let $x_+,x_-,h$ be the Chevalley basis of the simple Lie algebra
$sl_2$, that is
$$
[h,x_+]=2x_+,\;\; [h,x_-]=-2x_-,\;\; [x_+,x_-]=h.
$$
For any irreducible $sl_2$-module $V$ and $\alpha_1,\alpha_2\in{\bf
C}$, let $V(A)=V\otimes{\bf C}[t_1^{\pm 1},t_2^{\pm1}]$, we define
the action of $L$ on $V(A)$ as follows
$$\begin{array}{ll}
D(m_1{\bf e}_1+m_2{\bf b}_2)\cdot(v\otimes t_1^{n_1}t_2^{n_2})=
\!\!\!&(m_2(\alpha_1+n_1)-m_1(\alpha_2+n_2))v\otimes
t_1^{n_1+m_1}t_2^{n_2+m_2}
\\[5pt]&
+((m_2^2x_--m_1^2x_+m_1m_2h)\cdot v)\otimes
t_1^{n_1+m_1}t_2^{n_2+m_2}, \end{array}$$
$$
c_1.(v\otimes t_1^{n_1}t_2^{n_2})=c_2.(v\otimes
t_1^{n_1}t_2^{n_2})=0, \quad \forall v\in V.
$$
Then one can easily
see that $V(A)$ becomes a ${\bf Z}^2$-graded $L$-module. And all the
dimensions of the homogeneous subspaces of $V(A)$ are $\dim V$. The
following result was obtained in [15].

{\it {\bf Theorem 4.9} $([15])\quad (1)$ If $\dim V\geq3$ then
$V(A)$ is an irreducible $L$-module.

$(2)$ If $\dim V=2$ then $V(A)$ has a unique proper submodule when
$(\alpha_1,\alpha_2)\not\in{\bf Z}^2$, otherwise, $V(A)$ has two
proper submodule.

$(3)$ If $\dim V=1$ then $V(A)$ is an irreducible $L$-module when
$(\alpha_1,\alpha_2)\not\in{\bf Z}^2$, otherwise, $V(A)$ can be
decomposed into the direct sum of two irreducible $L$-module.}

\begin{center}
\item\section*{$\S 5$\quad \bf  ${\bf
Z}$-graded modules of the intermediate series over the Virasoro-like
algebra}
\end{center}

The main result in this section is the following.

{\it {\bf Theorem 5.1.}\quad There does not exist any nontrivial
${\bf Z}$-graded irreducible $L$-module of the intermediate series.}

{\bf Proof:}\quad If $V$ is a nontrivial ${\bf Z}$-graded
irreducible module of the intermediate series, then $c_1=c_2=0$ and
the action of $D({\bf b}_1)$ is nondegenerate. In fact, if it is
degenerate then there exists a vector $0\neq v_i\in V_i$ such that
$D({\bf b}_1)\cdot v_i=0$. Hence we have
$$
0=D({\bf b}_1)\cdot D({\bf b}_2)\cdot v_i=-\mbox{det}({\bf b}_1,{\bf
b}_2)D({\bf b}_1+{\bf b}_2)\cdot v_i,
$$
by the definition of  a ${\bf Z}$-graded $L$-module of the
intermediate series, which implies that $V$ is a generalize highest
weight module. Thus $V$ is not a uniformly bounded ${\bf Z}$-graded
module by Theorem 3.7 which is absurd.

Similarly, one can prove that the action of $D(-{\bf b}_1)$ is
nondegenerate. Therefore $\dim V_i=1$ for all $i\in{\bf Z}$. Thus we
can choose a basis $\{v_i\in V_i\,|\,i\in{\bf Z}\}$ of $V$ such that
$D({\bf b}_1)\cdot v_i=av_{i+1}$ and $D(-{\bf b}_1)\cdot
v_i=av_{i-1}$ for all $i\in{\bf Z}$, where $a\neq0$. Denote
$D((l{\bf b}_1+k{\bf b}_2)\cdot v_i=f(l,k,i)v_{l+i}$. We can deduce
$af(l,0,i\pm1)=af(l,0,i)$ since $D(l{\bf b}_1)D(\pm{\bf b}_1)\cdot
v_i=D(\pm{\bf b}_1)D(l{\bf b}_1)\cdot v_i$. Hence $f(l,0,i)$ are
independent of $i$ for all $l\neq0$. Set
$\varepsilon=\mbox{det}({\bf b}_1,{\bf b}_2)$ then
$\varepsilon\in\{\pm1\}$ since ${\bf b}_1,{\bf b}_2$ is a ${\bf
Z}$-basis of $\Gamma$.

For any $l,i,k\in{\bf Z}$ with $k\neq0$, we have
\begin{eqnarray*}&\!\!\!\!\!\!\!\!\!\!\!\!\!\!\!\!\!\!\!\!
 \varepsilon
kf(l,k,i)v_{l+i}\!\!\!&=\varepsilon kD(l{\bf b}_1+k{\bf b}_2)v_i
\\&\!\!\!\!\!\!\!\!\!\!\!\!\!\!\!\!\!\!\!\!&
=D((l-1){\bf b}_1+k{\bf b}_2)D({\bf b}_1)v_i-D({\bf b}_1)D((l-1){\bf
b}_1+k{\bf
b}_2)v_i\\&\!\!\!\!\!\!\!\!\!\!\!\!\!\!\!\!\!\!\!\!&=a(f(l-1,k,i+1)-f(l-1,k,i))v_{l+i},
\end{eqnarray*}
thus
$$\varepsilon f(l,k,i)=\frac{a}{k}(f(l-1,k,i+1)-f(l-1,k,i)).
\eqno{(5.1)}
$$
Similarly, we have
$$
\varepsilon kD((l-1){\bf b}_1+k{\bf b}_2)v_i=D(-{\bf b}_1)D(l{\bf
b}_1+k{\bf b}_2)v_i-D(l{\bf b}_1+k{\bf b}_2)D(-{\bf b}_1)v_i,
$$
which implies
$$
\varepsilon f(l-1,k,i)=\frac{a}{k}(f(l,k,i)-f(l,k,i-1)).
\eqno{(5.2)}
$$
Substituting (5.2) into (5.1), we have
$$
a^2f(l,k,i+1)-(2a^2+k^2)f(l,k,i)+a^2f(l,k,i-1)=0. \eqno{(5.3)}
$$
Set $x_k=\frac{2a^2+k^2+(4a^2k^2+k^4)^{\frac{1}{2}}}{2a^2}$. We can
choose $k\in{\bf N}$ such that $x_k\neq x_k^{-1}$ and $|x_k|>1$
since $\lim_{k\rightarrow +\infty}x_k=\infty$. Therefore, the
equation $a^2T^2-(2a^2+k^2)T+a^2=0$ has different roots $x_k$ and
$x_k^{-1}$, so we have
$$
f(l,k,i)=a(l,k)x_k^i+b(l,k)x_k^{-i},\quad \forall i\in{\bf Z},
\eqno{(5.4)}
$$
for some $a(l,k),b(l,k)\in{\bf C}$ by (5.3). Since $$\varepsilon
klD(l{\bf b}_1)v_i=D(l{\bf b}_1+k{\bf b}_2)D(-k{\bf
b}_2)v_i-D(-k{\bf b}_2)D(l{\bf b}_1+k{\bf b}_2)v_i,$$ we obtain
$$
\varepsilon klf(l,0,i)=f(l,k,i)(f(0,-k,i)-f(0,-k,l+i)).
$$
Thus, by (5.4) and the fact $x_k=x_{-k}$, we have
$$
\begin{array}{ll}\varepsilon
klf(l,0,i)=\!\!\!\!&a(l,k)a(0,-k)(1-x_k^l)x_k^{2i}+b(l,k)b(0,-k)(1-x_k^{-l})x_k^{-2i}
\\[5pt]&
+a(l,k)b(0,-k)(1-x_k^{-l})+b(l,k)a(0,-k)(1-x_k^l).
\end{array}\eqno{(5.5)}
$$ Since $|x_k|>1$, we see that if $a(l,k)a(0,-k)\neq0$ then
$$\begin{array}{ll}
\lim\limits_{i\rightarrow+\infty}(a(l,k)a(0,-k)(1-x_k^l)x_k^{2i}+b(l,k)b(0,-k)(1-x_k^{-l})x_k^{-2i}
\\[5pt]
+a(l,k)b(0,-k)(1-x_k^{-l})+b(l,k)a(0,-k)(1-x_k^l))=\infty,
\end{array}$$ which contradicts  (5.5) and the fact that $f(l,0,i)$ is
independent of $i$. Hence $a(l,k)a(0,-k)=0$. Similarly, we can
deduce $b(l,k)b(0,-k)=0$. If $a(0,-k)=b(0,-k)=0$ or
$a(l,k)=b(l,k)=0$ then $\varepsilon klf(l,0,i)\equiv 0$ by (5.5)
which is a contradiction. Hence without loss of generality, we may
assume that
$$
a(0,-k)=b(l,k)=0. \eqno{(5.6)}
$$
Thus, by (5.5) and (5.4), we have
$$
f(l,0,i)=\frac{a(l,k)b(0,-k)(1-x_k^{-l})}{\varepsilon lk},\quad
f(0,-k,i)=b(0,-k)x_k^{-i}, \quad f(l,k,i)=a(l,k)x_k^i. \eqno{(5.7)}
$$
Now we can substitute (5.7) into (5.1) to get
$$
\frac{a(l,k)}{a(l-1,k)}=a\varepsilon^{-1}k^{-1}(x_k-1).
$$
Thus $a(l,k)=a(0,k)(a\varepsilon^{-1}k^{-1}(x_k-1))^l$. From this
equation and (5.7), we can deduce that
$$f(l,k,i)=a(0,k)(a\varepsilon^{-1}k^{-1}(x_k-1))^lx_k^i,\;\;\forall
l\in{\bf Z},
 \eqno{(5.8)}
$$
and
$$
f(l,0,i)=\frac{a(0,k)b(0,-k)(1-x_k^{-l})(a\varepsilon^{-1}k^{-1}(x_k-1))^l}{\varepsilon
lk},\;\;\forall l\neq0. \eqno{(5.9)}
$$
Since $\varepsilon lkD(l{\bf b}_1+k{\bf b}_2)v_i=D(k{\bf
b}_2)D(l{\bf b}_1)v_i-D(l{\bf b}_1)D(k{\bf b}_2)v_i$, we have
$$
\varepsilon lkf(l,k,i)=a(0,k)x_k^i(x_k^l-1)f(l,0,i),
$$
by (5.8). Thus
$$
\frac{f(l,k,i)}{f(l,0,i)}=\frac{a(0,k)x_k^i(x_k^l-1)}{\varepsilon
lk}, \;\;\forall 0\neq l\in{\bf Z}.\eqno{(5.10)}
$$
Now by (5.8), (5.9) and (5.10), we have
$$
1=a(0,k)b(0,-k)l^{-2}k^{-2}(x_k^l-1)(1-x_k^{-l}),\;\;\forall 0\neq
l\in {\bf Z}.
$$
But $\lim_{l\rightarrow+\infty}(l^{-2}(x_k^l-1)(1-x_k^{-l}))
=\infty$ since $|x_k|>1$, this is a contradiction with  the above
equation. Thus there does not exist a nontrivial ${\bf Z}$-graded
irreducible $L$-module of the intermediate series.\hfill$\Box$

\end{document}